\documentclass[10pt]{amsart}
\usepackage[pagebackref]{hyperref}

\renewcommand{\d}{\mathrm{d}}

\newcommand{\ts}{\textstyle }

\newcommand{\cF}{{\mathcal F}}
\newcommand{\cG}{{\mathcal G}}

\newcommand{\cS}{{\mathcal S}}

\newcommand{\bbR}{{\mathbb R}}

\newcommand{\bbP}{{\mathbb P}}

\newcommand{\I}{\mathrm{I}}

\newcommand{\euso}{\operatorname{\mathfrak{so}}}
\newcommand{\euco}{\operatorname{\mathfrak{co}}}
\newcommand{\euspin}{\operatorname{\mathfrak{spin}}}
\newcommand{\eusp}{\operatorname{\mathfrak{sp}}}

\newcommand{\eugl}{\operatorname{\mathfrak{gl}}}
\newcommand{\eug}{\operatorname{\mathfrak{g}}}

\newcommand{\GL}{\operatorname{GL}}
\newcommand{\Spin}{\operatorname{Spin}}
\newcommand{\Symp}{\operatorname{Sp}}
\newcommand{\SO}{\operatorname{SO}}
\newcommand{\CO}{\operatorname{CO}}

\newcommand{\Gr}{\operatorname{Gr}}
\newcommand{\G}{\operatorname{G}}

\DeclareMathOperator{\tr}{tr}
\DeclareMathOperator{\Aut}{Aut}

\newcommand{\bth}{{\bar\theta}}
\newcommand{\bom}{{\bar\omega}}

\newcommand{\w}{{\mathchoice{\,{\scriptstyle\wedge}\,}{{\scriptstyle\wedge}}
      {{\scriptscriptstyle\wedge}}{{\scriptscriptstyle\wedge}}}}

\newcommand{\be}{\begin{equation}}
\newcommand{\ee}{\end{equation}}
\newcommand{\bpm}{\begin{pmatrix}}
\newcommand{\epm}{\end{pmatrix}}

\numberwithin{equation}{section}

\newtheorem{proposition}{Proposition}
\newtheorem{corollary}{Corollary}

\theoremstyle{remark}
\newtheorem{definition}{Definition}

\newtheorem{remark}{Remark}

\begin{document}

\author[R. Bryant]{Robert L. Bryant}
\address{Duke University Mathematics Department\\
         P.O. Box 90320\\
         Durham, NC 27708-0320}
\email{\href{mailto:bryant@math.duke.edu}{bryant@math.duke.edu}}
\urladdr{\href{http://www.math.duke.edu/~bryant}%
         {http://www.math.duke.edu/\lower3pt\hbox{\symbol{'176}}bryant}}

\title[Conformal geometry and $3$-plane fields]
      {Conformal geometry and $3$-plane fields \\
          on $6$-manifolds}

\date{December 29, 2005}

\begin{abstract}
The purpose of this note is to provide yet
another example of the link between certain
conformal geometries and ordinary differential
equations, along the lines of the examples
discussed by Nurowski~\cite{MR2157414}.

In this particular case, I consider the equivalence problem 
for~$3$-plane fields~$D\subset TM$ on a $6$-manifold~$M$
satisfying the nondegeneracy condition that $D+[D,D]=TM$.  

I give a solution of the equivalence problem for such~$D$
(as Tanaka has previously), showing that it defines 
a $\euso(4,3)$-valued Cartan connection 
on a principal right~$H$-bundle over~$M$ 
where~$H\subset\SO(4,3)$ is the subgroup
that stabilizes a null $3$-plane in~$\bbR^{4,3}$.
Along the way, I observe that there is associated
to each such~$D$ a canonical conformal structure of split type on~$M$, 
one that depends on two derivatives of the plane field~$D$.

I show how the primary curvature tensor of the
Cartan connection associated to the equivalence
problem for~$D$ can be interpreted
as the Weyl curvature of the associated conformal
structure and, moreover, show that the split conformal 
structures in dimension~$6$ that arise in this fashion
are exactly the ones whose $\euso(4,4)$-valued
Cartan connection admits a reduction 
to a~$\euspin(4,3)$-connection.  I also discuss 
how this case is analogous to features of Nurowski's examples.
\end{abstract}



\thanks{
Thanks to Duke University for its support via a research grant, 
to the National Science Foundation 
for its support via grants DMS-9870164 and DMS-0103884, 
and to the Clay Mathematics Institute for its support 
in the period January-May 2002, during which a
portion of this manuscript was written.\hfill\break
\hspace*{\parindent} 
This is Version~$2.0$. 
The most recent version can be found at \texttt{arXiv:math.DG/0511110}.
}


\maketitle

\setcounter{tocdepth}{2}
\tableofcontents

\section{Introduction}\label{sec:intro}

In~\cite{MR2157414}, Nurowski considers several
different equivalence problems for classes of differential
equations and shows how each one leads to a natural 
conformal structure (of indefinite type) of an 
appropriate configuration space and that this
conformal structure suffices to encode 
the original equivalence problem.

Perhaps the most striking of these examples is the
one based on \'E. Cartan's famous `five-variables'
paper~\cite{Cartan1910}, in which Cartan solves
the equivalence problem for $2$-plane fields of
maximal growth vector~$(2,3,5)$ on $5$-manifolds.   
Such $2$-plane fields are now said to be `of Cartan type'
in honor of Cartan's pioneering work. 

In that paper, Cartan shows that, 
given such a $2$-plane field~$D\subset TM$
where~$M$ has dimension~$5$, one can define
what is now called a Cartan connection over~$M$
that solves the equivalence problem.  Specifically,
let~$\G_2'\subset\SO(4,3)$ be the noncompact
exceptional simple group of dimension~$14$.  The group~$\G_2'$
acts transitively on the set~$Q_{3,2}\simeq S^3\times S^2$
of null lines in~$\bbR^{4,3}$.  Let~$H\subset \G_2'$
be the subgroup of codimension~$5$ that fixes a null
line in~$Q_{3,2}$.  Then Cartan shows how canonically 
to associate to~$D$ a principal right~$H$-bundle~$\pi:P\to M$
and a $\eug_2'$-valued $1$-form~$\gamma$ on~$P$ 
such that each (possibly locally defined)
diffeomorphism~$\phi:M\to M$ that preserves~$D$ lifts
canonically to an $H$-bundle diffeomorphism~$\hat\phi:P\to P$
that fixes~$\gamma$.  Cartan shows, further, that 
part of the curvature of~$\gamma$ can be interpreted
as a section~$\cG$ of the bundle~$S^4(D_1^*)$, 
where~$D_1 = D + [D,D]$ is the rank~$3$ first derived bundle of~$D$.  
He also shows that the necessary and sufficient condition for `flatness',
i.e., equivalence of~$D$ with the ${\rm G}_2'$-invariant
$2$-plane field on~$Q_{3,2}$, is that this section of~$S^4(D_1^*)$
should vanish identically.  In fact, he proves the stronger fact 
that~$\cG$ vanishes if and only if the `restricted'
curvature, i.e., the reduced section of~$S^4(D^*)$, which
Cartan denotes as~$\cF$, vanishes.
(Recall that, since the inclusion~$D\to D_1$ is an injection, 
the dual restriction map~$S^4(D_1^*)\to S^4(D^*)$ is a surjection.)

Of course, $\G_2'$ preserves a conformal structure
of split type on~$Q_{3,2}$.  What Nurowski shows is that,
for general~$D$ of Cartan type, there is associated 
a natural conformal structure of split type on~$M$, 
generalizing the case of~$Q_{3,2}$. 
He also shows that Cartan's tensor~$\cG$
is simply the Weyl curvature of this associated conformal structure.  

In this note, I point out a similar result for $3$-plane fields
on $6$-manifolds~$D\subset TM$ that satisfy the generic condition
that~$D+[D,D]=TM$.  

In~\S\ref{sec:equiv_prob}, 
I work out the equivalence problem for such $3$-plane fields.
Of course, following the work of Cartan, this is just a calculation.
Moreover, Tanaka~\cite{MR0221418,MR0533089} has explained how 
to solve this problem (and many more like it), so this aspect 
of the article is not at all new.%
\footnote{In fact, \S\ref{sec:equiv_prob} and~\S\ref{sec:fund_tensor}
are based on calculations that I did in my 1979 thesis~\cite{Bryant1979}, 
when I was ignorant of Tanaka's work.  These sections were
actually written for a series of lectures that I gave in~2002
on the method of equivalence but never published.}

One thing that is, perhaps, 
new, and is motivated by comparison with Nurowski's work, is the
observation, made in Proposition~\ref{eq:splitconformalexist}, that
there exists a canonical conformal structure of split type on~$M$ 
associated to such a $3$-plane field.  This conformal structure
depends on two derivatives of the defining equations of the $3$-plane
field, as is evidenced by the fact that it is first defined in
terms of the second-order frame bundle, as derived in the 
course of the equivalence problem.

The result of the equivalence problem calculation is that,
if~$H\subset\SO(4,3)$ is the stabilizer subgroup of a null $3$-plane
in~$\bbR^{4,3}$, then the plane field~$D\subset TM$ defines a
principal right~$H$-bundle~$B_3\to M$ and a $\euso(4,3)$-valued
Cartan connection $1$-form~$\gamma$ on~$B_3$ such that every diffeomorphism
$\phi:M\to M$ that preserves the plane field~$D$ induces in a 
canonical way a lifted $H$-bundle automorphism~$\hat\phi:B_3\to B_3$
that preserves the Cartan connection~$\gamma$.  Moreover,
every $H$-bundle map~$\varphi:B_3\to B_3$ that preserves~$\gamma$
is of the form~$\varphi = \hat\phi$ for a unique diffeomorphism~$\phi:M\to M$
that preserves~$D$.

I show that the fundamental curvature tensor of~$\gamma$, 
which I denote by~$\cS$, can be regarded as a section 
of the rank~$27$ Shur-irreducible bundle 
\be
\bigl(S^2(D)\otimes S^2(D^*)\bigr)_0\otimes\Lambda^3(D^*).
\ee
This fundamental curvature tensor is the analog of Cartan's reduced curvature, 
i.e., in his case, the section of~$S^4(D^*)$ 
(his `binary quartic form'~$\cF$) 
rather than of~$S^4(D_1^*)$ (his `ternary quartic form'~$\cG$).
Correspondingly, in this case, there is, in fact, an extended
curvature tensor~$\cS^+$ that has a canonical reduction 
to~$\cS$, but I do not write it out explicitly here.

I show that the vanishing of~$\cS$ is the necessary and
sufficient condition that~$D$ be locally equivalent to the `flat
example', i.e., the $3$-plane field on~$\SO(4,3)/H$ that is 
preserved by the action of~$\SO(4,3)$.  (In particular, I show
that the vanishing of~$\cS$ implies that of~$\cS^+$.)

Finally, I show that the tensor~$\cS^+$ 
is simply the Weyl curvature of the conformal structure 
on~$M$ associated to~$D$,
exactly as Nurowski shows in Cartan's case.

\section{The equivalence problem}\label{sec:equiv_prob}

\subsection{Maximal non-integrability}\label{ssec:max_non-int}
Let~$M$ be a smooth $6$-manifold and let~$D\subset TM$
be a smooth $3$-plane field with the property that
the set~$D+[D,D]$ is equal to~$TM$ and has constant rank.
In other words, every point~$x\in M$ has
a neighborhood~$U$ on which there exist vector fields~$X_1$,
$X_2$, $X_3$ that are sections of~$D$ over~$U$, are everywhere 
linearly independent on~$U$, and have the property
that the six vector fields
\be\label{eq:bracket_span}
X_1\,,\ X_2\,,\ X_3\,,\ [X_2,X_3]\,,\ [X_3,X_1]\,,\ [X_1,X_2]
\ee 
are everywhere linearly indepdendent on~$U$.  
Thus,~$D$ is `maximally nonintegrable'.  

A dual formulation
of this maximal non-integrability condition is that there 
exist $1$-forms $\theta_1$, $\theta_2$, and~$\theta_3$ on~$U$ 
so that each~$\theta_i$ annihilates all of the vectors in~$D$ and
so that~$\d\theta_1$, $\d\theta_2$, and $\d\theta_3$ are
linearly independent modulo~$\theta_1$, $\theta_2$, and~$\theta_3$
everywhere on~$U$.

\subsection{$1$-adaptation}\label{ssec:1_adaptation}
A coframing~$\eta:TU\to\bbR^6$ on an open set~$U\subset M$ 
of the form
\be\label{eq:coframing_eta}
\eta = \begin{pmatrix}
\bth_1\\ \bth_2\\ \bth_3\\ 
\bom^1\\ \bom^2\\ \bom^3
\end{pmatrix}
\ee
will be said to be \emph{$1$-adapted}
to~$D$ if each of the $\bth_i$ annihilate the vectors in~$D$
and if the equations
\be\label{eq:1-adapt}
\left.
\begin{aligned}
\d\bth_1 &\equiv 2\,\bom^2\w\bom^3\\
\d\bth_2 &\equiv 2\,\bom^3\w\bom^1\\
\d\bth_3 &\equiv 2\,\bom^1\w\bom^2\\
\end{aligned}
\ \right\} \mod \bth_1,\bth_2,\bth_3
\ee
hold on~$U$.

The coframings $1$-adapted to~$D$ are the local sections of 
a $G_1$-structure~$B_1\to M$, where~$G_1\subset\GL(6,\bbR)$ 
is the group of matrices of the form
\be\label{eq:G_1-definition}
\begin{pmatrix}
\det(A)\, {}^t\!A^{-1} & 0 \\
AB & A 
\end{pmatrix}
\ee
where~$A$ lies in~$\GL(3,\bbR)$ 
and~$B$ is an arbitrary $3$-by-$3$ matrix.

I will denote the entries of the tautological $\bbR^6$-valued $1$-form
on~$B_1$ as~$\theta_i$ and~$\omega^i$, as in equation~\eqref{eq:coframing_eta}.
By construction, there exists on~$B_1$ a pseudo-connection of the form
\be\label{eq:pseudo-conn}
\begin{pmatrix}
\tr(\alpha)\,\I_3 - {}^t\alpha  & 0 \\
\beta & \alpha
\end{pmatrix}
\ee
where~$\alpha=(\alpha^i_j)$ and~$\beta=(\beta^{ij})$ take
values in $3$-by-$3$-matrices, so that equations of the form%
\footnote{Here, as henceforth in this article, the summation convention
is to be assumed.}
\be\label{eq:zeroth_str_eqs}
\begin{aligned}
\d\theta_i &= -\alpha^k_k\w\theta_i + \alpha^j_i\w\theta_j 
              + \epsilon_{ijk}\,\omega^j\w\omega^k\\
\d\omega^i &= -\beta^{ij}\w\theta_j -\alpha^i_j\w\omega^j 
              + P^{il}\epsilon_{ljk}\,\omega^j\w\omega^k\\
\end{aligned}
\ee
hold, where~$P^{il}$ are some functions on~$B_0$ and where~$\epsilon_{ijk}$
is totally skewsymmetric in its indices and satisfies~$\epsilon_{123}=1$.

\subsection{$2$-adaptation and a conformal structure}
\label{ssec:2_adaptation}
Now, expanding out~$\d(\d\theta_i)=0$ and reducing the result modulo
$\theta_1$, $\theta_2$, and~$\theta_3$ yields the relations~$P^{il}
=P^{li}$.  One now finds that the six equations~$P^{il}=0$ define
a sub-bundle~$B_2\subset B_1$ that is a~$G_2$-structure on~$M$,
where~$G_2\subset G_1$ is the subgroup%
\footnote{Of course, the reader will not confuse~$G_2$ 
with~${\rm G}_2$, the simple group of dimension~$14$.} 
consisting of those
matrices of the form~\eqref{eq:G_1-definition} in which~$B$ is
skewsymmetric, i.e.,~${}^t\!B = - B$.  A coframing~$\eta$ that
is a section of~$B_2$ will be said to be \emph{$2$-adapted} to~$D$.

\begin{proposition}\label{eq:splitconformalexist}
There exists a unique pseudo-conformal structure of split type on~$M$ 
such that a nondegenerate quadratic form~$g$ on~$M$ 
represents this conformal structure if and only if its
pullback to~$B_2$ is a multiple of the quadratic 
form~$\theta_i\circ\omega^i$.
\end{proposition}

\begin{proof}
Note the evident fact that~$G_2$ 
is a subgroup of the group~$\CO(3,3)\subset\GL(6,\bbR)$
consisting of the invertible matrices~$h$ that satisfy
\be
{}^th\bpm 0_3&I_3\\I_3 &0_3\epm h
= \bigl|\det(h)\bigr|^{1/3} \bpm 0_3&I_3\\I_3 &0_3\epm.
\ee
The proposition now follows since~$B_2$ is a $G_2$-structure on~$M$.
\end{proof}

\begin{remark}[Order of the conformal structure]
Note that, because the bundle~$B_2$ is constructed out of
two derivatives of the plane field~$D$, the conformal structure
depends on two derivatives of the plane field~$D$.
\end{remark}

\begin{remark}[A weighted quadratic form]
In fact, one can get a well-defined tensor on~$M$ 
out of this construction:  Let~$\eta:U\to B_2$ be 
a $2$-adapted coframing on a domain~$U\subset M$ 
and write~$\eta$ in the form~\eqref{eq:coframing_eta}.  
Let~$X_i$ be the sections of~$D$ over~$U$ 
that satisfy~$\bom^i(X_j) = \delta^i_j$.  
Then the tensor
\be
\hat g = \bth_i{\circ}\bom^i\otimes (X_1\w X_2\w X_3)
\ee
is a well-defined section of~$S^2(T^*M)\otimes\Lambda^3(D)$
that depends on two derivatives of~$D$.  Clearly, $\hat g$
determines the canonical conformal structure.
\end{remark}

Pulling the pseudo-connection forms back to~$B_2$ and writing
$\beta^{ij} = \epsilon^{ikj}\beta_k + \tau^{ij}$ where~$\tau^{ij}
=\tau^{ji}\equiv0\mod\{\theta,\omega\}$, 
the structure equations take the form
\be\label{eq:first_str_eqs0}
\begin{aligned}
\d\theta_i &= -\alpha^k_k\w\theta_i + \alpha^j_i\w\theta_j 
              + \epsilon_{ijk}\,\omega^j\w\omega^k\,,\\
\d\omega^i &= -\epsilon^{ikj}\,\beta_k\w\theta_j -\alpha^i_j\w\omega^j 
              + \tau^{ij}\w\theta_k\,.
\end{aligned}
\ee

Setting
\be\label{eq:A_defn1}
A^i_j = \d\alpha^i_j + \alpha^i_k\w\alpha^k_j + 2\,\omega^i\w\beta_j
\ee
and expanding the identity~$\d(\d\theta_i)=0$ now yields
\be\label{eq:ddtheta=0}
0 = -A^k_k\w\theta_i + A^j_i\w\theta_j 
\ee
from which it follows, in particular, that~$A^i_j\equiv0\mod\{\theta\}$.
Using this congruence to expand the identity~$\d(\d\omega^i)=0$
and then reducing modulo~$\{\theta\}$ yields
\be\label{eq:ddomega=0modtheta}
0\equiv \tau^{ij}\w \epsilon_{jkl}\,\omega^k\w\omega^l \mod\{\theta\}.
\ee
It follows that there exist functions~$T^{ij}_k=T^{ji}_k$ that
satisfy~$T^{ij}_i=0$ and
\be\label{eq:tau_mod_theta}
\tau^{ij} \equiv T^{ij}_k\,\omega^k
\mod\{\theta\}.
\ee
Now, by a replacement of the form~$\alpha^i_j\mapsto \alpha^i_j+p^{ij}_k\,
\theta^k$, one can retain the first equations of~\eqref{eq:first_str_eqs0}
(this imposes $9$ linear equations on the $27$ functions~$p^{ij}_k$)
while simultaneously reducing the functions~$T^{ij}_k$ to zero
(this imposes $15$ further linear equations on the $27$ functions~$p^{ij}_k$
and these are independent from the first~$9$).

Thus, there exist pseudo-connection forms~$\alpha^i_j$ and~$\beta_j$
on~$B_2$ so that the equations
\be\label{eq:first_str_eqs1}
\begin{aligned}
\d\theta_i &= -\alpha^k_k\w\theta_i + \alpha^j_i\w\theta_j 
              + \epsilon_{ijk}\,\omega^j\w\omega^k\\
\d\omega^i &= -\epsilon^{ikj}\,\beta_k\w\theta_j -\alpha^i_j\w\omega^j 
              + \epsilon^{ljk}T^i_l\,\theta_j\w\theta_k
\end{aligned}
\ee
hold for some functions~$T^i_j$ on~$B_2$.  
However, again, by linear algebra, there exists a unique
replacement of the form~$\beta_i\mapsto\beta_i+p^j_i\theta_j$
for which~$T^i_j=0$.  Thus, there exist pseudo-connection 
forms~$\alpha^i_j$ and~$\beta_j$ on~$B_2$ so that the equations
\be\label{eq:first_str_eqs2}
\begin{aligned}
\d\theta_i &= -\alpha^k_k\w\theta_i + \alpha^j_i\w\theta_j 
              + \epsilon_{ijk}\,\omega^j\w\omega^k\\
\d\omega^i &= -\epsilon^{ikj}\,\beta_k\w\theta_j -\alpha^i_j\w\omega^j
\end{aligned}
\ee
hold.  The pseudo-connection forms are not uniquely determined
by these equations; one can perform the replacements
\be\label{eq:pdo_conn_var}
\begin{aligned}
\alpha^i_j &\longmapsto \alpha^i_j + \delta^i_j\,t^k\,\theta_k - t^i\,\theta_j\\
\beta_i &\longmapsto \beta_i +\epsilon_{ijk}\,t^j\,\omega^k
\end{aligned}
\ee
for any functions~$t^1,t^2,t^3$ without affecting~\eqref{eq:first_str_eqs2}.
(Of course, this corresponds to the fact that the first 
prolongation~$\eug_2^{(1)}$ of the algebra~$\eug_2\subset\eugl(6,\bbR)$ 
has dimension~$3$.)

\subsection{Prolongation and the third-order bundle}\label{ssec:prolongation}
Let~$B_3\to B_2$ be the $\bbR^3$-bundle over~$B_2$ whose fibers are the
point pseudo-connections for which equations~\eqref{eq:first_str_eqs2}
hold.  Then equations
\be\label{eq:second_str_eqs0}
\begin{aligned}
\d\theta_i &= -\alpha^k_k\w\theta_i + \alpha^j_i\w\theta_j 
              + \epsilon_{ijk}\,\omega^j\w\omega^k\\
\d\omega^i &= -\epsilon^{ikj}\,\beta_k\w\theta_j -\alpha^i_j\w\omega^j
\end{aligned}
\ee
hold on~$B_3$, where now the forms~$\theta$, $\omega$, $\alpha$, and~$\beta$
are tautologically defined (and, hence, canonical).  Set
\be\label{eq:AB_defn1}
\begin{aligned}
A^i_j &= \d\alpha^i_j + \alpha^i_k\w\alpha^k_j + 2\,\omega^i\w\beta_j\\
B_i &= \d\beta_i - \alpha^j_i\w\beta_j.
\end{aligned}
\ee
The exterior derivatives of the equations~\eqref{eq:second_str_eqs0} 
can now be expressed as
\be\label{eq:dsecond_str_eqs0}
\begin{aligned}
0 &= -A^k_k\w\theta_i + A^j_i\w\theta_j \\
0 &= -\epsilon^{ikj}\,B_k\w\theta_j - A^i_j\w\omega^j.
\end{aligned}
\ee
The first equation of~\eqref{eq:dsecond_str_eqs0} implies, 
in particular, that~$A^i_j\equiv 0 \mod\{\theta\}$, so there exist
$1$-forms~$\pi^{ik}_j$ (not unique) such that~$A^i_j=\pi^{ik}_j\w\theta_k$.
Substituting this relation into the second set of equations
of~\eqref{eq:dsecond_str_eqs0} then yields
\be\label{eq:Bpi_reln}
0 = -\epsilon^{ikj}\,B_k\w\theta_j - \pi^{ij}_k\w\theta_j\w\omega^k,
\ee
which, in turn, implies
\be\label{eq:Bpi_reln_modtheta}
0 \equiv -\epsilon^{ikj}\,B_k + \pi^{ij}_k\w\omega^k \mod\{\theta\}.
\ee
In particular, it follows 
that~$\pi^{ij}_k+\pi^{ji}_k\equiv0\mod\{\theta,\omega\}$,
so that one can write $\pi^{ij}_k=\epsilon^{ijl}\pi_{lk}+\sigma^{ij}_k$
where~$\sigma^{ij}_k=\sigma^{ji}_k\equiv0\mod\{\theta,\omega\}$.  One
can further write~$\pi_{ij} = -\epsilon_{ijk}\,\tau^k+\sigma_{ij}$ where
$\sigma_{ij}=\sigma_{ji}$.  This leads to the formula
\be
A^i_j = \pi^{ik}_j\w\theta_k 
= \delta^i_j\,\tau^k\w\theta_k - \tau^i\w\theta_j 
+\epsilon^{ikl}\sigma_{jl}\w\theta_k +\sigma^{ik}_j\w\theta_k.
\ee
Substituting this into the first set of equations in~\eqref{eq:dsecond_str_eqs0}
and using the fact that $\sigma^{ij}_k\equiv0\mod\{\theta,\omega\}$ shows
that the $3$-forms~$\Sigma_j =\sigma_{jl}\w\epsilon^{ikl}\,\theta_k\w\theta_i$
are cubic expressions in the $1$-forms~$\theta_i$ and~$\omega^k$.  
In particular, it follows that~$\sigma_{jl}\equiv0\mod\{\theta,\omega\}$.
Consequently, the $2$-forms~$A^i_j$ can be written in the form
\be\label{eq:Acurv0}
A^i_j =  \delta^i_j\,\tau^k\w\theta_k - \tau^i\w\theta_j 
         + R^i_{jk}\epsilon^{klm}\,\theta_l\w\theta_m
	      + S^{ik}_{jl}\,\theta_k\w\omega^l
\ee
for some $1$-forms~$\tau^i$ and functions~$R^i_{jk}$ and~$S^{ik}_{jl}$.
Comparing this with~equation~\eqref{eq:pdo_conn_var}, one sees that
the 1-forms~$\tau^i$ are the components of a pseudo-connection 
for the bundle~$B_3\to B_2$.  Of course, these~$\tau^i$ are not uniquely
determined by the formulae~\eqref{eq:Acurv0}.  

\subsubsection{Normalizing~$\tau$.}\label{sssec:taunormalization}
The $\tau_i$ will be made unique by imposing the appropriate 
linear equations on the functions~$R$ and~$S$ as follows:  
First, consider the trace of~\eqref{eq:Acurv0}:
\be\label{eq:Acurv0traced}
A^i_i =  2\,\tau^k\w\theta_k 
         + R^i_{ik}\epsilon^{klm}\,\theta_l\w\theta_m
	      + S^{ik}_{il}\,\theta_k\w\omega^l\,.
\ee
By adding linear combinations of the $\omega^i$ and $\theta_j$ 
to the $\tau^k$, one can arrange that
\be\label{eq:RStrace_conditions}
R^i_{ik} = S^{ik}_{il} = 0.
\ee 
In other words, $A^i_i = 2\,\tau^k\w\theta_k$.  

The conditions~\eqref{eq:RStrace_conditions} 
still do not determine the~$\tau^k$ completely.  
However, they do determine the~$\tau^k$ up to a replacement 
of the form~$\tau^k\mapsto\tau^k+p^{kl}\,\theta_l$ 
where~$p^{kl}=p^{lk}$.  

Substituting these normalized formulae into the first set
of equations in~\eqref{eq:dsecond_str_eqs0} 
yields the relations
\be
0 =\left(R^j_{ik}\epsilon^{klm}\,\theta_l\w\theta_m
	      + S^{jk}_{il}\,\theta_k\w\omega^l\right)\w\theta_j\,,
\ee
which are equivalent to the equations
\be\label{eq:RS_1stBianchi}
R^j_{ij} = S^{jk}_{il}-S^{kj}_{il} = 0.
\ee 

This suggests a closer inspection of the functions~$R^i_{jk}$.  
Consider the $\GL(3,\bbR)$-invariant decomposition
\be\label{eq:Rdecomp0}
R^i_{jk} = S^i_{jk} + \epsilon_{ljk}\,S^{il} + \epsilon_{ljk}\epsilon^{ilp}S_p\,,
\ee
where~$S^i_{jk}=S^i_{kj}$ and~$S^{ij}=S^{ji}$.  
The trace condition~$R^i_{ij}=0$ and identity~$R^i_{ji}=0$ 
now combine to show that~$S^i_{ij}=S_j=0$, 
so the decomposition of~$R$ simplifies to
\be\label{eq:Rdecomp1}
R^i_{jk} = S^i_{jk} + \epsilon_{ljk}\,S^{il}\,,
\ee
where~$S^i_{jk}=S^i_{kj}$ and~$S^{ij}=S^{ji}$.

One can now finally complete the normalization of the~$\tau^k$ 
by requiring, in addition to~\eqref{eq:RStrace_conditions}, that~$S^{ij}=0$.
Thus, the~$\tau^k$ are made unique by requiring them to be chosen so that
\be\label{eq:Acurv1}
A^i_j =  \delta^i_j\,\tau^k\w\theta_k - \tau^i\w\theta_j 
         + S^i_{jk}\epsilon^{klm}\,\theta_l\w\theta_m
	      + S^{ik}_{jl}\,\theta_k\w\omega^l
\ee
holds, where the coefficients are required to satisfy the normalizations
\be\label{eq:Snormalizations}
S^i_{jk}=S^i_{kj}\,,\ S^i_{ik}=0\,,\ S^{ik}_{ij}=0
\ee

Thus, the forms~$\theta_i$, $\omega^j$, $\alpha^i_j$, $\beta_i$, and~$\tau^j$
define a canonical coframing on~$B_3$ and every diffeomorphism of~$M$
that preserves the $3$-plane field~$D$ lifts to a unique diffeomorphism
of~$B_3$ that fixes the forms in this coframing.  Thus, these constitute
the solution of the equivalence problem in the sense of Cartan.  

\subsection{Bianchi identities}\label{ssec:further_identities}
Substituting equation~\eqref{eq:Acurv1} into the second
set of equations of~\eqref{eq:dsecond_str_eqs0}, 
yields the relations 
\be
\epsilon^{ikj}\,B_k\w\theta_j =
    - \left(\delta^i_j\,\tau^k\w\theta_k - \tau^i\w\theta_j 
         + S^i_{jk}\epsilon^{klm}\,\theta_l\w\theta_m
	      + S^{ik}_{jl}\,\theta_k\w\omega^l\right)\w\omega^j.
\ee
It follows that~$S^{ik}_{jl}=S^{ik}_{lj}$ 
and that one has relations of the form
\be\label{Bcurv0}
B_i = \epsilon_{ijk}\,\tau^j\w\omega^k - 2 S^j_{ik}\,\theta_j\w\omega^k
        +\epsilon^{jkl}S_{ij}\,\theta_k\w\theta_l\,,
\ee
where~$S_{ij}=S_{ji}$.  

To summarize the results so far:  There are structure equations
\be\label{eq:second_str_eqs1}
\begin{aligned}
\d\theta_i &= -\alpha^k_k\w\theta_i + \alpha^j_i\w\theta_j 
              + \epsilon_{ijk}\,\omega^j\w\omega^k\\
\d\omega^i &= -\epsilon^{ikj}\,\beta_k\w\theta_j -\alpha^i_j\w\omega^j\\
\d\alpha^i_j  &= -\alpha^i_k\w\alpha^k_j - 2\,\omega^i\w\beta_j +  
            \delta^i_j\,\tau^k\w\theta_k - \tau^i\w\theta_j \\
        &\qquad\qquad + S^i_{jk}\epsilon^{klm}\,\theta_l\w\theta_m
	      + S^{ik}_{jl}\,\theta_k\w\omega^l\\ 
\d\beta_i &= \alpha^j_i\w\beta_j + \epsilon_{ijk}\,\tau^j\w\omega^k 
        - 2 S^j_{ik}\,\theta_j\w\omega^k +\epsilon^{jkl}S_{ij}\,\theta_k\w\theta_l
\end{aligned}
\ee
where the functions~$S$ satisfy the trace and symmetry relations
\be\label{eq:Srelations}
S_{ij}=S_{ji}\,,\ \ S^i_{jk}=S^i_{kj}\,,\ S^i_{ik}=0\,,
\ \ S^{ik}_{jl}=S^{ki}_{jl}=S^{ik}_{lj}\,,\ S^{ik}_{ij}=0.
\ee

Tracing the formula for~$\d\alpha^i_j$ yields
\be\label{eq:Acurvtraced}
\d\alpha^i_i = 2\,\tau^k\w\theta_k - 2\,\omega^k\w\beta_k
\ee
and, since the left-hand side is closed, 
taking the exterior derivative of both sides yields
\be
0 = 2\,T^k \w \theta_k
\ee
where
\be
T^i = \d\tau^i - \alpha^k_k\w\tau^i + \alpha^i_j\w\tau^j
      + \epsilon^{ijk}\,\beta_j\w\beta_k
	  + \epsilon^{ijl}S_{lk}\,\theta_j\w\omega^k\,.
\ee
Thus, there exist $1$-forms~$\tau^{ij}=\tau^{ji}$ so that
\be\label{eq:dtau0}
\d\tau^i  = \alpha^k_k\w\tau^i - \alpha^i_j\w\tau^j
      - \epsilon^{ijk}\,\beta_j\w\beta_k
	  - \epsilon^{ijl}S_{lk}\,\theta_j\w\omega^k
	  - \tau^{ij}\w\theta_j\,.
\ee 
These $1$-forms~$\tau^{ij}$ are not unique, 
but are unique up to a replacement 
of the form~$\tau^{ij}\mapsto\tau^{ij}+ p^{ijk}\,\theta_k$
for some functions~$p^{ijk}$ satisfying 
the symmetry conditions~$p^{ijk}=p^{jik}=p^{ikj}$.

Define $1$-forms~$\sigma_{ij}$, $\sigma^i_{jk}$, 
and $\sigma^{ik}_{jl}$ by the equations
\be\label{eq:sigma_defn}
\begin{aligned}
\d S^{ik}_{jl} &= \sigma^{ik}_{jl} + S^{ik}_{jl}\,\alpha^m_m
   -S^{mk}_{jl}\alpha^i_m -S^{im}_{jl}\alpha^k_m
   +S^{ik}_{ml}\alpha^m_j +S^{ik}_{jm}\alpha^m_l\\
   &\quad +{\ts\frac23}\left(5\delta^k_mS^i_{jl}+5\delta^i_mS^k_{jl}
                       -\delta^k_jS^i_{ml}- \delta^i_jS^k_{ml}
					   -\delta^k_lS^i_{jm}- \delta^i_lS^k_{jm}
                       \right)\omega^m \\
\d S^i_{jk} 
&= \sigma^i_{jk} + S^i_{jk}\,\alpha^m_m - S^m_{jk}\,\alpha^i_m
    + S^i_{mk}\,\alpha^m_j+S^i_{jm}\,\alpha^m_k+{\ts\frac12}S^{il}_{jk}\,\beta_l\\
&\quad - {\ts\frac12}\left(4\delta^i_lS_{jk}
            -\delta^i_kS_{jl}-\delta^i_jS_{lk}\right)\omega^l
			 \\
\d S_{ij}  &= \sigma_{ij} + S_{jk}\,\alpha^m_m 
+ S_{mj}\,\alpha^m_i+S_{im}\,\alpha^m_j-2S^{m}_{ij}\,\beta_m\,. \\
\end{aligned}
\ee
Then the~$\sigma$s satisfy the same symmetry and trace conditions
as the corresponding~$S$s and, moreover, the identies~$\d(\d\alpha^i_j)=0$
and~$\d(\d\beta_i)=0$ become the relations
\be\label{eq:dda_ddb}
\begin{aligned}
0  &=  - \tau^{im}\w\theta_m\w\theta_j 
       + \epsilon^{klm}\,\sigma^i_{jk}\w\theta_l\w\theta_m
	      + \sigma^{ik}_{jl}\w\theta_k\w\omega^l\\ 
0 &= \epsilon_{ijk}\,\tau^{jl}\w\theta_l\w\omega^k 
        - 2 \sigma^j_{ik}\w\theta_j\w\omega^k
		+\epsilon^{jkl}\sigma_{ij}\w\theta_k\w\theta_l
\end{aligned}
\ee
These relations imply
\be
\tau^{ij}\equiv\sigma_{ij}\equiv\sigma^i_{jk}\equiv\sigma^{ik}_{jl}
\equiv 0 \mod\{\theta,\omega\}.
\ee
(If~$X$ is a vector field on~$B_3$ that satisfies~$\theta_i(X)=\omega^j(X)=0$,
then the above equations imply
\be\label{eq:dda_ddb_lhkX}
\begin{aligned}
0  &=  - \tau^{im}(X)\,\theta_m\w\theta_j 
       + \epsilon^{klm}\,\sigma^i_{jk}(X)\,\theta_l\w\theta_m
	      + \sigma^{ik}_{jl}(X)\,\theta_k\w\omega^l\\ 
0 &= \epsilon_{ijk}\,\tau^{jl}(X)\,\theta_l\w\omega^k 
        - 2 \sigma^j_{ik}(X)\,\theta_j\w\omega^k
		+\epsilon^{jkl}\sigma_{ij}(X)\,\theta_k\w\theta_l\,,
\end{aligned}
\ee
which implies $\tau^{ij}(X)=\sigma_{ij}(X)=\sigma^i_{jk}(X)=\sigma^{ik}_{jl}(X)=0$.
Hence the conclusion.) 

It follows that there are expansions
\be\label{eq:sigma_expansions}
\begin{aligned}
  \tau^{ij}_{  } & = T^{ijm}_{ }\,\theta_m + T^{ij\phantom{m}}_{m}\,\omega^m\\
\sigma^{  }_{ij} & = T^{m\phantom{ik}}_{ij}\,\theta_m + T^{ }_{ijm}\,\omega^m\\
\sigma^{i }_{jk} & = T^{im\phantom{k}}_{jk}\,\theta_m + T^{i}_{jkm}\,\omega^m\\
\sigma^{ik}_{jl} & = T^{ikm}_{jl}\,\theta_m + T^{ik}_{jlm}\,\omega^m\\
\end{aligned}
\ee
and that~$\tau^{ij}$ can be made unique 
by requiring that the full symmetrization of~$T^{ijm}$ vanish, 
i.e., that $T^{ijk}+T^{jki}+T^{kij}=0$, 
so assume that this has been done.

The relations~\eqref{eq:dda_ddb} 
can now be expressed as the following identities:
\be\label{eq:T_identities}
\begin{aligned}
T^{ik}_{jlm} & = T^{ik}_{jml}\\
T^i_{jkm} &= {\ts\frac14}\left(\epsilon_{jpq}\,T^{ipq}_{km}+
                    \epsilon_{kpq}\,T^{ipq}_{jm}\right)\\
T^{ij}_m &= -{\ts\frac12}\,T^{ijk}_{mk}\\
T^{ijm} &= {\ts\frac23}\left(\epsilon^{ilm}\,T^{jk}_{lk}+
                    \epsilon^{jlm}\,T^{ik}_{lk}\right)\\
T^m_{im} &= 0\\
T_{ijm} &= 3\epsilon_{klm}\,T^{kl}_{ij}-2\epsilon_{kli}\,T^{kl}_{mj}
           -2\epsilon_{klj}\,T^{kl}_{mi}+\epsilon_{ilm}\,T^{lk}_{jk}
		   +\epsilon_{jlm}\,T^{lk}_{ik}\,.
\end{aligned}
\ee

\section{The fundamental tensor and flatness}\label{sec:fund_tensor}
The expansions~\eqref{eq:sigma_expansions} taken with the definition
of~$\sigma^{ik}_{jl}$ in~\eqref{eq:sigma_defn} show that the
functions~$S^{ij}_{kl}$ are constant on the fibers of~$B_3\to B_2$
and hence can be regarded as functions on~$B_2$.  In fact, because
\be\label{eq:dSikjk-variance}
\d S^{ik}_{jl} \equiv S^{ik}_{jl}\,\alpha^m_m
   -S^{mk}_{jl}\alpha^i_m -S^{im}_{jl}\alpha^k_m
   +S^{ik}_{ml}\alpha^m_j +S^{ik}_{jm}\alpha^m_l \mod \{\theta,\omega\},
\ee
it follows that the~$S^{ik}_{jl}$ can be regarded as
the components of a section of the 
bundle~$S^2(D)\otimes S^2(D^*)\otimes\Lambda^3(D^*)$ 
that takes values in the (irreducible) Shur representation
subbundle~$\bigl(S^2(D)\otimes S^2(D^*)\bigr)_0\otimes\Lambda^3(D^*)$,
which has rank~$27$ (the subscript~$0$ denotes the kernel of the 
natural mapping~$S^2(D)\otimes S^2(D^*)\to D\otimes D^*$ 
that is defined by contraction).

Specifically, if~$\eta = (\bth_i,\bom^j)$ is a $2$-adapted coframing
on some domain~$U\subset M$, set~${\bar S}^{ik}_{jl} = \eta^*S^{ik}_{jl}$
and consider the expression
\be
\cS(\eta) = {\bar S}^{ik}_{jl}\,\,{\bar X}_i{\circ}{\bar X}_k
              \otimes \bom^j{\circ}\,\bom^l
              \otimes \left(\bom^1{\w}\bom^2{\w}\bom^3\right)
\ee
as a section of~$S^2(D)\otimes S^2(D^*)\otimes\Lambda^3(D^*)$ over~$U$,
where~${\bar X}_i$ are the sections of~$D$ over~$U$ that are dual
to~$\bom^i$, i.e., so that~$\bom^i({\bar X}_j)=\delta^i_j$.  
Then equation~\eqref{eq:dSikjk-variance} implies that~$\cS(\eta)$ 
is independent of the choice of $1$-adapted coframing~$\eta$ 
and hence is the restriction to~$U$ 
of a globally defined section~$\cS$ 
that depends only on~$D$.

\begin{definition}[The fundamental tensor]
The tensor~$\cS$ will be referred to 
as the \emph{fundamental tensor} of~$D$.
\end{definition} 

The following vanishing result is the analog for nondegenerate 
$3$-plane fields in dimension~$6$ 
of Cartan's characterization in~\cite[\S VII]{Cartan1910}
of the `flat' $2$-plane fields of Cartan type in dimension~$5$.

\begin{proposition}\label{prop:Svanishes}
Suppose that~$\cS$ vanishes identically.
Then the following hold: First,~$S^i_{jk}$ and~$S_{ij}$
vanish identically.  Second, the structure 
equations simplify to
\be\label{eq:flat_str_eqs}
\begin{aligned}
\d\theta_i &= -\alpha^k_k\w\theta_i + \alpha^j_i\w\theta_j 
              + \epsilon_{ijk}\,\omega^j\w\omega^k\\
\d\omega^i &= -\epsilon^{ikj}\,\beta_k\w\theta_j -\alpha^i_j\w\omega^j\\
\d\alpha^i_j  &= -\alpha^i_k\w\alpha^k_j - 2\,\omega^i\w\beta_j +  
            \delta^i_j\,\tau^k\w\theta_k - \tau^i\w\theta_j \\
\d\beta_i &= \phantom{-}\alpha^j_i\w\beta_j + \epsilon_{ijk}\,\tau^j\w\omega^k\\
\d\tau^i  &= \phantom{-}\alpha^k_k\w\tau^i - \alpha^i_j\w\tau^j
      - \epsilon^{ijk}\,\beta_j\w\beta_k\,.
\end{aligned}
\ee
Third, for any $1$-connected open~$U\subset M$, the Lie algebra 
of vector fields on~$U$ whose {\upshape(}local{\upshape)} flows 
preserve~$D$ is isomorphic to the Lie algebra of~$\SO(4,3)$.
Fourth, any point of~$M$ is the center of a coordinate system~$\bigl(U,
(x^j,y_i)\bigr)$ in which the plane field~$D$ is annihilated by the three
$1$-forms~$\bth_i = \d y_i + \epsilon_{ijk}\,x^j\,\d x^k$.
\end{proposition}

\begin{proof}
First, note that, by the first equation of~\eqref{eq:sigma_defn}, 
the vanishing of the functions~$S^{ik}_{jl}$ implies that
\be
\begin{aligned}
\sigma^{ik}_{jl}
&= -{\ts\frac23}\left(5\delta^k_mS^i_{jl}+5\delta^i_mS^k_{jl}
                       -\delta^k_jS^i_{ml}- \delta^i_jS^k_{ml}
                       -\delta^k_lS^i_{jm}- \delta^i_lS^k_{jm}
                       \right)\omega^m\\
&= T^{ikm}_{jl}\,\theta_m + T^{ik}_{jlm}\,\omega^m\,
\end{aligned}
\ee
which, in turn, implies both~$T^{ikm}_{jl}=0$ and
\be
T^{ik}_{jlm}
= -{\ts\frac23}\left(5\delta^k_mS^i_{jl}+5\delta^i_mS^k_{jl}
                       -\delta^k_jS^i_{ml}- \delta^i_jS^k_{ml}
                       -\delta^k_lS^i_{jm}- \delta^i_lS^k_{jm}\right).
\ee
However, by the first equation of~\eqref{eq:T_identities}, $T^{ik}_{jlm}$
is fully symmetric in its lower indices, which implies
\be
S^i_{jk} = 0.
\ee
Using this, by the second equation of~\eqref{eq:sigma_defn}
and by~\eqref{eq:sigma_expansions}, one has
that
\be
\sigma^i_{jk} = {\ts\frac12}\left(4\delta^i_lS_{jk}
                    -\delta^i_kS_{jl}-\delta^i_jS_{lk}\right)\omega^l
               = T^{im\phantom{k}}_{jk}\,\theta_m + T^{i}_{jkl}\,\omega^l\,
\ee
which implies that~$T^{im\phantom{k}}_{jk} = 0$ and
\be
T^{i}_{jkl} = {\ts\frac12}\left(4\delta^i_lS_{jk}
                    -\delta^i_kS_{jl}-\delta^i_jS_{lk}\right).
\ee
Now, however, the second equation of~\eqref{eq:T_identities} coupled
with~$T^{ikm}_{jl}=0$ (which was derived above) show that~$T^{i}_{jkl} = 0$
which, in turn, now implies~$S_{ij}=0$. 

Next, since~\eqref{eq:sigma_defn} now implies that~$\sigma_{ij}=0$,
it follows from~\eqref{eq:sigma_expansions}, that
\be
T^{m\phantom{ik}}_{ij} =  T^{ }_{ijm} = 0.
\ee
A final appeal to~\eqref{eq:T_identities} then shows that
\be
T^{ijm}_{ } =  T^{ij\phantom{m}}_{m} = 0,
\ee
i.e., that~$\tau^{ij}=0$.  Consequently, the structure equations
simplify to~\eqref{eq:flat_str_eqs}, as claimed.

Now, the exterior derivatives of the equations~\eqref{eq:flat_str_eqs}
are identities, so it follows that these are the left-invariant
forms on a Lie group of dimension~$21$.  An examination of the 
weights associated to the (maximal) torus dual to the diagonal~$\alpha$s
shows that the Lie algebra is~$\euso(4,3)$, thus
verifying the third claim.  

One can also see this directly by noting that the
equations~\eqref{eq:flat_str_eqs} are equivalent 
to~$\d\gamma = -\gamma\w\gamma$, where
\be
\gamma
= \bpm \vspace{3pt}
-\alpha^1_1&-\alpha^2_1&-\alpha^3_1&2\beta_1&0&-\tau_3&\tau_2\\ \vspace{3pt}
-\alpha^1_2&-\alpha^2_2&-\alpha^3_2&2\beta_2&\tau_3&0&-\tau_1\\ \vspace{3pt}
-\alpha^1_3&-\alpha^2_3&-\alpha^3_3&2\beta_3&-\tau_2&\tau_1&0\\ \vspace{3pt}
  \omega^1&\omega^2&\omega^3&0&-\beta_1&-\beta_2&-\beta_3\\ \vspace{3pt}
  0&\theta_3&-\theta_2&-2\omega^1&\alpha^1_1&\alpha^1_2&\alpha^1_3\\ \vspace{3pt}
  -\theta_3&0&\theta_1&-2\omega^2&\alpha^2_1&\alpha^2_2&\alpha^2_3\\ \vspace{3pt}
  \theta_2&-\theta_1&0&-2\omega^3&\alpha^3_1&\alpha^3_2&\alpha^3_3
  \epm .
\ee
Obviously, $\gamma$ takes values 
in the Lie algebra~$\euso(4,3)\subset\eugl(7,\bbR)$, which is the
space of matrices~$a$ that satisfy~$Qa + {}^taQ = 0$,
where
\be
Q = \bpm 0_{3\times3} & 0_{3\times1} & I_{3\times3}\\
         0_{1\times3} & 2 & 0_{1\times3}\\ 
         I_{3\times3} & 0_{3\times1} & 0_{3\times3}\epm
\ee
is a symmetric matrix of type~$(4,3)$.

Finally, the system~$\alpha^i_j=\beta_j=\tau^i=0$ is a Frobenius
system, and a leaf of this system in~$B_3$ defines a (local)~$2$-adapted 
coframing~$\eta$ on an open set~$U\subset M$ that satisfies
\be
\begin{aligned}
\d\bom^j &= 0\,,\\
\d\bth_{i} &= \epsilon_{ijk}\,\bom^j\w\bom^k\,.
\end{aligned}
\ee
Consequently, assuming that~$U$ is simply connected, there exist
functions~$x^j$ on~$U$ such that~$\bom^j = \d x^j$ and there exist
functions~$y_i$ on~$U$ such that
\be
\d y_i=\bth_i-\epsilon_{ijk}\,x^j\,\d x^k.
\ee
These provide the desired local coordinates.
\end{proof}

\begin{corollary}[Maximal symmetry] 
The Lie group~$\Aut(M,D)$ has dimension at most~$21$
and this upper limit is reached only when~$D$ is locally
equivalent to the $3$-plane field on~$\bbR^6$ defined by
the equations
\be
\d y_i+\epsilon_{ijk}\,x^j\,\d x^k = 0.
\ee
\end{corollary}

\begin{proof}
The Lie group~$\Aut(M,D)$ is embedded into the group of 
diffeomorphisms of~$B_3$ that preserve the 
coframing defined by~$\theta_i$,~$\omega^j$,~$\alpha^j_i$,~$\beta_j$,
and~$\tau^i$.  This group can only have dimension~$21$ if all of
the functions~$S^{ik}_{jl}$ are constant.  However, these 
functions cannot be constant on the fibers of~$B_3\to M$ unless
they vanish.  Now apply Proposition~\ref{prop:Svanishes}.
\end{proof}

\begin{remark}[The homogeneous model]
Note that the proof of Proposition~\ref{prop:Svanishes} identifies
the homogeneous model for the `flat' case:  Let~$M^6\subset\Gr(3,\bbR^{4,3})$ 
be the space of isotropic (i.e., null) $3$-planes in the split signature
inner product space~$\bbR^{4,3}$.  The group~${\rm O}(4,3)$ acts
transitively on this $6$-manifold and preserves a nondegenerate $3$-plane
field on it.  By Proposition~\ref{prop:Svanishes}, the identity
component of~${\rm O}(4,3)$ is the identity component of 
the group of automorphisms of this $3$-plane field.
\end{remark}

\begin{remark}[Irregular $D$-curves]
Note that the Cartan system of the $1$-form~$\theta_3$ on~$B_3$
is the Pfaffian system~$J$ spanned by~$\theta_1$,
$\theta_2$, $\theta_3$, $\omega^1$, $\omega^2$, $\alpha^1_3$, and
$\alpha^2_3$.  Consequently, this Pfaffian system is Frobenius
(as can be directly verified by a glance at the structure equations)
and hence there is a submersion~$\nu:B_3\to N^7$ for some (not 
necessarily Hausdorf) $7$-manifold~$N^7$ such that the fibers of~$\nu$
are the leaves of~$J$.  

The points of~$N^7$ represent the \emph{irregular $D$-curves} 
in~$M^6$, as defined in~\cite{MR1240644}.%
\footnote{Actually, in \cite{MR1240644}, these curves are
called `non-regular', but I now prefer the more standard English
term `irregular'.}
Specifically, a leaf of the system~$J$ projects to~$M$ 
as a submersion onto a curve in~$M$ and, in this way, 
one sees that each~$J$-leaf represents a curve in~$M$.  

This $7$-parameter family of curves has the property that
exactly one curve of the family passes through a given point in~$M$
with a given tangent direction in~$D$.  

Note that, in the homogeneous model, $N^7$ is simply the space
of null (i.e., isotropic) $2$-planes in~$\bbR^{4,3}$.  
Each such $2$-plane lies in a $1$-parameter family of null $3$-planes
and this gives the interpretation of such $2$-planes as curves in~$M$.
In fact, given a null $2$-plane~$E\subset\bbR^{4,3}$, the restriction
of the quadratic form to the $5$-plane~$E^\perp\subset\bbR^{4,3}$
has kernel equal to~$E$ and hence descends to a nondegenerate
form (of type~$(2,1)$ on~$E^\perp/E\simeq\bbR^{2,1}$.  The null $3$-planes
that contain~$E$ are in $1$-to-$1$ correspondence with the null lines
in~$E^\perp/E$, a space which is known to be $1$-dimensional and,
in fact, naturally isomorphic to~$\bbR\bbP^1$.

Similarly, in the general case, each of the irregular $D$-curves 
inherits a natural projective structure.  In fact, on a leaf of~$J$, 
one has $\d\omega^3 = -\alpha^3_3\w\omega^3$,
$\d\alpha^3_3 = 2\,\beta_3\w \omega^3$, and~$\d\beta_3 = \alpha^3_3\w\beta_3$,
so that~$\omega^3$ is a differential on the corresponding $D$-curve 
that is well-defined up to a projective change of parameter.
\end{remark}

\begin{remark}[An extended tensor]
The reader cannot have helped but notice that 
equations~\eqref{eq:sigma_defn} actually imply that
$\cS$ is the reduction of an extended tensor~$\cS^+$
of rank~$48 = 27{+}15{+}6$ 
that uses all of the components~$S^{ij}_{kl}$, $S^i_{kl}$, and~$S_{kl}$.
This extended tensor will play a role in the next section, but it is
not worthwhile to write it out explicitly here.  
Instead, I will just note that~$\cS^+$ 
takes values in a certain rank~$48$ subbundle of
the bundle~$S^2(\eug_2)\otimes\Lambda^3(D^*)$, 
where~$\eug_2\subset\eugl(6,\bbR)$ is the Lie algebra of the subgroup~$G_2$
defined at the beginning of~\S\ref{ssec:2_adaptation}.

For a comparison with Nurowski's examples, 
see~Remark~\ref{rem:Nurowski5dimlcase}.
\end{remark}

\section{The fundamental tensor and Weyl curvature}\label{sec:Weylinterp}
Consider the~$1$-form~$\hat\gamma$ with values in~$\euso(4,4)\subset\eugl(8,\bbR)$
defined on~$B_3$ by the formula
\be\def\phm{\phantom{-}}
\hat\gamma
 = 
\bpm\vspace{3pt}
-\phi &\phm\beta_1&\phm\beta_2&\phm\beta_3
&\phm\tau^1&\phm\tau^2&\phm\tau^3&0\\ \vspace{3pt}
\phm\omega^1&\alpha^1_1{-}\phi&\alpha^1_2&\alpha^1_3&0&-\beta_3&\phm\beta_2
&\phm\tau^1\\ \vspace{3pt}
\phm\omega^2&\alpha^2_1&\alpha^2_2{-}\phi&\alpha^2_3&\phm\beta_3&0&-\beta_1
&\phm\tau^2\\ \vspace{3pt}
\phm\omega^3&\alpha^3_1&\alpha^3_2&\alpha^3_3{-}\phi&-\beta_2&\phm\beta_1&0
&\phm\tau^3\\ \vspace{3pt}
\phm\theta_1&0&\phm\omega^3&-\omega^2&\phi{-}\alpha^1_1&-\alpha^2_1&-\alpha^3_1
&\phm\beta_1\\ \vspace{3pt}
\phm\theta_2&-\omega^3&0&\phm\omega^1&-\alpha^1_2&\phi{-}\alpha^2_2&-\alpha^3_2
&\phm\beta_2\\ \vspace{3pt}
\phm\theta_3&\phm\omega^2&-\omega^1&0&-\alpha^1_3&-\alpha^2_3&\phi{-}\alpha^3_3
&\phm\beta_3\\  \vspace{3pt}
\phm0&\phm\theta_1&\phm\theta_2&\phm\theta_3
&\phm\omega^1&\phm\omega^2&\phm\omega^3&\phm\phi
\epm
\ee
where~$\phi = \frac12(\alpha^i_i)$.

In the flat case, this $1$-form satisfies~$\d\hat\gamma
=-\hat\gamma\w\hat\gamma$.  In particular, it follows that~$\hat\gamma$ takes
values in a Lie algebra~$\eug\subset\euso(4,4)$ that is isomorphic to~$\euso(4,3)$.
The corresponding subgroup of~$\SO(4,4)$ is isomorphic to~$\Spin(4,3)$.
Thus, I will denote the algebra~$\eug$ by~$\euspin(4,3)$ 
and call the corresponding subgroup~$\Spin(4,3)$.

It now follows from~\eqref{eq:second_str_eqs1} that, 
if~$P\to M$ is the Cartan structure bundle 
associated to the canonical conformal structure 
with~$\euso(4,4)$-valued connection form~$\Gamma$ 
and fiber isomorphic to the parabolic subgroup~$H\subset\SO(4,4)$ 
that is the stabilizer of a null line in~$\bbR^{4,4}$,
then there exists a bundle embedding~$\iota:B_3\to P$ such that
\be
\hat\gamma = \iota^*(\Gamma).
\ee 
In particular, the structure equations~\eqref{eq:second_str_eqs1}
show that the functions~$S^{ik}_{jl}$, $S^i_{jl}$, $S_{jl}$ are
the components of the Weyl curvature of the conformal structure
in this reduction.  Thus, one has the following result:

\begin{proposition}\label{prop:Weyl}
The Weyl tensor of the conformal structure associated to~$D$
is the extended tensor~$\cS^+$.  In particular, the
associated conformal structure is conformally flat if and only
if the plane field~$D$ is locally equivalent to the 
flat example.\qed
\end{proposition}

\begin{remark}[An algebraic characterization of the Weyl tensor]
Recall that, for a split-conformal manifold of dimension~$6$, 
the Weyl tensor takes values in a bundle associated to an 
irreducible, $84$-dimensional representation space~$W$ of~$\euco(3,3) 
= \bbR\oplus\euso(3,3)$ that can be described as follows:  
Use the `exceptional isomorphism' $A_3 = D_3$ to regard~$\euco(3,3)$
as~$\eugl(4,\bbR)$ and let~$V$ be the standard representation of
dimension~$4$ of~$\eugl(4,\bbR)$.  Then it is not difficult to establish
the isomorphism of representations
\be
W = \bigl(S^2(V){\otimes}S^2(V^*)\bigr)_0\otimes\bigl(\Lambda^4(V)\bigr)^{-1/2}
\ee 
where~$\bigl(S^2(V){\otimes}S^2(V^*)\bigr)_0\subset S^2(V){\otimes}S^2(V^*)$
is the kernel of the natural (and surjective) contraction mapping
\be
S^2(V)\otimes S^2(V^*) \longrightarrow V\otimes V^*.
\ee

Now, if~$\xi\subset V^*$ is a hyperplane, one can define the subspace
\be
\bigl(S^2(V)\otimes S^2(\xi)\bigr)_0\subset\bigl(S^2(V) \otimes S^2(V^*)\bigr)_0
\ee
to be the kernel of the natural (and surjective) contraction mapping
\be
S^2(V)\otimes S^2(\xi) \longrightarrow V\otimes \xi.
\ee
The dimension of this space is~$48$ and it is a representation
space of the~$12$-dimensional subgroup~$G_\xi\subset\GL(V)$ 
that preserves the hyperplane~$\xi$.  Under the isomorphism 
(actually, a double cover)~$\GL(V)\to\CO(3,3)$, 
the subgroup~$G_\xi$ goes to the subgroup~$G_2\subset\CO(3,3)$.

Now, the Weyl curvature function of the conformal structure
pulls back to~$B_3$ to take values in the $48$-dimensional subspace
\be
W_\xi = \bigl(S^2(V){\otimes}S^2(\xi)\bigr)_0\otimes\bigl(\Lambda^4(V)\bigr)^{-1/2},
\ee
This subspace is characterized as the kernel of the contraction
\be
C_e:W\to \bigl(S^2(V)\otimes V^*)_0\otimes\bigl(\Lambda^4(V)\bigr)^{-1/2}
\ee
where~$e\subset V$ is a nonzero vector annihilated by~$\xi$.  
  
The group~$G_\xi$ preserves the filtration
\be
S^2(\xi^\perp)\otimes S^2(\xi)
\subset\bigl(\xi^\perp{\circ}V\otimes S^2(\xi)\bigr)_0
\subset\bigl(S^2(V)\otimes S^2(\xi)\bigr)_0
\ee
whose graded pieces have dimensions~$6$, $15$, and~$27$.  
This filtration corresponds to the representation of the Weyl curvature 
by the components~$S_{jk}$, $S^i_{jk}$ and~$S^{il}_{jk}$, which
are the components of the tensor~$\cS^+$.  
In particular, the top associated graded piece
\be
\frac{\bigl(S^2(V)\otimes S^2(\xi)\bigr)_0}
     {\bigl(\xi^\perp{\circ}V\otimes S^2(\xi)\bigr)_0}
\simeq \bigl(S^2(V/\xi^\perp)\otimes S^2(\xi)\bigr)_0
\ee
of dimension~$27$ gives the associated bundle 
in which the tensor~$\cS$ takes values.

Thus, the algebraic characterization of the Weyl tensors 
that arise from conformal structures associated 
to nondegenerate $3$-plane fields on $6$-manifolds 
is that there should exist a nonzero conformal half-spinor~$s$ 
whose contraction with the Weyl tensor should vanish.
This nonvanishing half-spinor field then defines the structure reduction 
that locates~$B_3$ as a subbundle of the conformal Cartan connection bundle.
\end{remark}

\begin{remark}[Analogy with the $5$-dimensional case]
\label{rem:Nurowski5dimlcase}
By Nurowski's calculations in~\cite{MR2157414},
Proposition~\ref{prop:Weyl} has a direct parallel 
in the case of Cartan-type $2$-plane fields in dimension~$5$.

In that case, Nurowski shows that Cartan's ternary quartic form~$\cG$
can be interpreted as the Weyl curvature of the conformal structure
associated to the $2$-plane field.

In fact, the parallel is even more striking when one looks at the
algebraic characterization of the Weyl curvature in that case.
There, Cartan's structure bundle~$\pi:P\to M^5$ and $\eug_2'$-valued
connection form~$\gamma$ are embedded via an equivariant inclusion
$\iota:P\to P^+$ into the conformal structure bundle~$\pi:P^+\to M^5$ 
with $\euso(4,3)$-valued connection form~$\Gamma$.  This corresponds
to the inclusion~${\rm G}_2'\subset\SO(4,3)$ 
and the fiber group~$H\subset {\rm G}_2'$ of the bundle~$P$ 
is the intersection of ${\rm G}_2'$ with the subgroup~$H^+\subset\SO(4,3)$ 
that consists of those elements that fix a given null line in~$\bbR^{4,3}$.

The group~$H^+$ has a natural homomorphism onto~$\CO(3,2)$ 
and the image of~$H$ under this natural homomorphism 
is the $7$-dimensional subgroup~$K\subset\CO(3,2)$ 
that fixes a null $2$-plane.

Now, the Weyl curvature of a conformal structure of split type 
in dimension~$5$ is an irreducible, $35$-dimensional representation~$W$
of~$\CO(3,2)$.  Using the exceptional isomorphism~$\euco(3,2)
\simeq\bbR\oplus\eusp(2,\bbR)$ and letting~$V$ denote the irreducible
$4$-dimensional representation of~$\bbR^*{\cdot}\Symp(2,\bbR)$, 
then~$W$ is isomorphic to~$S^4(V^*)$.  

Also, the subgroup of~$\SO(3,2)$ that fixes a null $2$-plane
corresponds, in~$\Symp(2,\bbR)$ to the $7$-dimensional
subgroup that fixes a line in~$V$, or, equivalently, 
a $3$-plane~$\xi\subset V^*$.  

Tracing through this isomorphism
and comparing it with the calculations of Nurowski, one sees
that the Weyl curvature of the conformal structure associated to
a Cartan-type $2$-plane field in dimension~$5$ takes values in
a subspace of the form~$S^4(\xi)\subset S^4(V^*)\simeq W$ for 
an appropriately chosen~$\xi$.  Moreover, since the subgroup of
$\Symp(2,\bbR)$ that fixes~$\xi$ preserves the 
subspace~$\xi^\perp\subset\xi$, it follows that this subgroup
preserves a filtration of~$S^4(\xi)$ based on the number of
factors of~$\xi^\perp$ that appear in the quartic.  This
filtration has graded pieces of degrees~$1$, $2$, $3$, $4$, and~$5$,
which corresponds exactly to the filtration of Cartan's tensor~$\cG$
into the components that he labels~$E$, $D_i$, $C_i$, $B_i$, and~$A_i$
(with $1$, $2$, $3$, $4$, and~$5$ components, respectively.)
\end{remark}

\bibliographystyle{hamsplain}

\begin{thebibliography}{10}

\bibitem{Bryant1979}
R. Bryant,
\emph{Some aspects of the local and global theory of Pfaffian systems},
PhD Thesis, University of North Carolina at Chapel Hill, 1979.

\bibitem{MR1240644}
R. Bryant and L. Hsu, 
\emph{Rigidity of integral curves of rank~$2$ distributions}, 
Inventiones Mathematicae \textbf{114} (1993), 435--461.
MR1240644 (94j:58003) 

\bibitem{Cartan1910}
\'E.~Cartan,  
\emph{Les syst\`emes de Pfaff a cinq variables 
           et les \'equations aux d\'eriv\'ees partielles du second ordre},
Ann.\ Sc.\ Norm.\ Sup.\ \textbf{27} (1910), 109--192. 


\bibitem{MR2157414}
P. Nurowski, 
\emph{Differential equations and conformal structures},
J.\ Geom.\ Phys. \textbf{55} (2005), 19--49.  math.DG/0406400

\bibitem{MR0221418}
N. Tanaka, 
\emph{On generalized graded Lie algebras 
         and geometric structures.~I}, 
J.\ Math.\ Soc.\ Japan \textbf{19} (1967),  215--254.
MR0221418 (36 \#4470) 

\bibitem{MR0533089} 
\bysame, 
\emph{On the equivalence problems
       associated with simple graded Lie algebras}, 
Hokkaido Math.\ J.\ \textbf{8} (1979),  23--84.
MR0533089 (80h:53034)

\end{thebibliography}

\providecommand{\bysame}{\leavevmode\hbox to3em{\hrulefill}\thinspace}

\end{document}